\documentclass[12pt]{amsart} 
\usepackage[mathscr]{eucal}
\usepackage{amsmath,amsfonts}

\parskip=\smallskipamount

\newtheorem{theorem}{Theorem}[section]

\newtheorem{example}[theorem]{Example}

\newtheorem{algorithm}[theorem]{Algorithm}

\makeatletter
\@addtoreset{equation}{section}
\makeatother

 \newcommand{\cA}{{\mathcal A}}

 \newcommand{\cD}{{\mathcal D}}

 \newcommand{\cH}{{\mathcal H}}
 
 \newcommand{\cL}{{\mathcal L}}
 \newcommand{\cM}{{\mathcal M}}

 \newdimen\expt
 \def\boxit#1{\setbox0\hbox{$\displaystyle{#1}$}
       \hbox{\lower.4\expt
  \hbox{\lower3\expt\hbox{\lower\dp0
       \hbox{\vbox{\hrule height.4\expt
  \hbox{\vrule width.4\expt\hskip3\expt
       \vbox{\vskip3\expt\box0\vskip2\expt}%
  \hskip3\expt\vrule width.4\expt}\hrule height.4\expt}}}}}}
 
\begin{document}
 \pagestyle{plain}

 \bigskip

 \title 
 {Lattice structures for quantum channels} 
 \author{T. Constantinescu}  

 \address{Department of Mathematics \\
   University of Texas at Dallas \\
   Box 830688, Richardson, TX 75083-0688, U. S. A.}
 \email{\tt tiberiu@utdallas.edu}

 \begin{abstract}
 We suggest that a certain one-to-one parametrization
 of completely positive maps on the matrix algebra $\cM_n$
 might be useful in the study of quantum channels. This is 
 illustrated in the case of binary quantum channels. While the algorithm 
 is quite intricate, it admits a simple, lattice structure
 representation. 

 \end{abstract}

 \maketitle

 \section{Introduction}

 Some recent papers deal with the analysis of 
 the completely positive, trace-preserving linear maps
 on the matrix algebra $\cM_n$, \cite{FA}, \cite{RSW}. 
 The analysis is quite complete in the case $n=2$, as it
 can be seen in the paper \cite{RSW}.
 The purpose of this paper is to introduce an algorithm 
 that tests the complete positivity of a linear map
 on $\cM_n$, for any $n\geq 2$. This appears as a sort of Schur-Cohn test
 and it allows the introduction of certain lattice structures
 associated to completely positive linear maps.
 The algorithm is applied to $\cM_2$ and the result is compared
 with the analysis in \cite{RSW}.
 Since our algorithm produces a "free" parametrization 
 of the completely positive maps on $\cM_n$, it is 
 nonlinear in nature and other applications 
 in order to check its usefulness remain to be investigated.

 \section{Completely positive maps on $\cM_n$}  
 Let $\cM _n$ denote the algebra of complex $n\times n$
 matrices. A linear map $\Phi :\cA\rightarrow \cL(\cH)$
 from a $C^*$-algebra into the set $\cL(\cH)$ of all 
 bounded linear operators on the Hilbert space $\cH$ 
 is called {\it completely positive}
 if for every positive integer $n$, the map
 $$\Phi \otimes I_{\cM _n}:\cA\otimes \cM_n\rightarrow
 \cL(\cH)\otimes \cM_n$$
 is positivity preserving. 
 By Stinespring Theorem, \cite{St}, \cite{Pa}, 
 any such map is the compression of a $*$-homomorphism. For linear
 completely positive maps on $\cM_n$, this implies a somewhat more
 explicit representation of the form:
 \begin{equation}\label{kraus}
 \Phi (X)=\sum _{j}A_j^*XA_j,
 \end{equation}
 where $\{A_j\}$ is a finite set of elements in $\cM_n$,
 \cite{Kr}, \cite{Ch}, \cite{Pa}. The representation \eqref{kraus}
 is not unique and another characterization of linear
 completely positive maps on $\cM_n$ is also useful. Thus, by
 a result of Choi \cite{Ch}, \cite{Pa},
 the linear map $\Phi :\cM_n\rightarrow \cM_n$
 is completely positive if and only if the matrix
 \begin{equation}\label{baza}
 S=S_{\Phi }=\left[
 \begin{array}{ccc}
 \Phi (E_{11}) & \ldots & \Phi (E_{1n}) \\
  & & \\
 \vdots & \ddots & \\
  & & \\
 \Phi (E_{n1}) &  & \Phi (E_{nn})
 \end{array}\right]
 \end{equation}
 is positive (semi-definite), where $\{E_{kj}\}_{k,j=1}^n$
 are the standard matrix units of $\cM_n$, that is, 
 $E_{kj}$ is $1$ in the $(k,j)$-th entry and $0$ elsewhere. 
 We notice that if $X=\left[X_{kj}\right]_{k,j=1}^n$,
 then $Y=\left[Y_{kj}\right]_{k,j=1}^n=\Phi (X)$ is given by the relations
 \begin{equation}\label{filtru}
 Y_{kj}=\sum _{l,m}\Phi (E_{lm})_{kj}X_{lm},
 \end{equation} 
 and the correspondence \eqref{baza} between the completely
 positive maps on $\cM_n$ and the $n^2\times n^2$ positive matrices
 is one-to-one and affine (see \cite{FA} for details).

 Of special interest in quantum information are those linear completely
 positive maps that preserve the trace. Such maps are usually called
 {\it quantum channels}, \cite{BS}.
 The adjoint $\hat{\Phi }$ of a linear map $\Phi $ on $\cM_n$ 
 is defined with respect to the Hilbert structure on $\cM_n$
 given by the Hilbert-Schmidt inner product
 (linear in the first variable),
 $\langle A,B\rangle =TrAB^*$, $A,B\in \cM_n$,
 where $B^*$ denotes the usual adjoint of the operator $B$.
 It follows that $\Phi $ is trace-preserving if and only if
 $\hat{\Phi }$ is unital ($\hat{\Phi }(I)=I$).

 We will use some standard notation associated to contractions 
 on Hilbert spaces. Thus, let $\cL(\cH_1,\cH_2)$ denote the set
 of all bounded linear maps operators from the Hilbert space 
 $\cH_1$ into the Hilbert space $\cH_2$.
 The operator $T\in \cL(\cH_1,\cH_2)$ is called a {\it contraction}
 if $\|T\|\leq 1$. The {\it defect operator} of $T$ is 
 $D_T=(I-T^*T)^{1/2}$ and $\cD_T$ denotes the closure of the range
 of $D_T$. To any contraction $T\in \cL(\cH_1,\cH_2)$
 one associates the unitary operator $U(T):\cH_1\oplus \cD_{T^*}
 \rightarrow \cH_2\oplus \cD_{T}$ by the formula:
 \begin{equation}\label{rotele}
 U(T)=\left[
 \begin{array}{cc}
 T & D_{T^*} \\
 D_T & -T^* 
 \end{array}
 \right].
 \end{equation}

 \section{Lattice structures}
 Let $\Phi $ be a linear completely positive map on $\cM_n$.
 The matrix $S=S_{\phi}$ given by \eqref{baza}
 is positive, and by Theorem~1.5.3 in 
 \cite{Co}, there exists a uniquely determined family
 $\Gamma =\{\Gamma _{kj}\mid 1\leq k\leq j\leq n^2\}$
 of complex numbers with the following properties.
 Thus, 
 $$S_{kk}=\Gamma _{kk}, \quad 1\leq k\leq n^2,$$
 and for $1\leq k<j\leq n^2$,
 $\Gamma _{kj}\in \cL(\cD_{\Gamma _{k+1,j}},\cD_{\Gamma ^*_{k,j-1}})$
 are contractions such that 
 \begin{equation}\label{algoritm}
 S_{kj}=\Gamma ^{1/2}_{kk}(R_{k,j-1}U_{k+1,j-1}C_{k+1,j}+
 D_{\Gamma ^*_{k,k+1}}\ldots 
 D_{\Gamma ^*_{k,j-1}}\Gamma _{kj}
 D_{\Gamma _{k+1,j}}\ldots
 D_{\Gamma _{j-1,j}})\Gamma ^{1/2}_{jj}.
 \end{equation}
 We use the convention that 
 $\cD_{\Gamma _{kk}}$ is just (the closure of) the range 
 of $\Gamma _{kk}$. This algorithm is valid in higher dimensions as
 well, that is the entries of $S$ could be bounded operators and 
 then the parameters $\Gamma _{kj}$ would be also operators.
 The notation used in \eqref{algoritm} is quite involved but easy 
 to explain. Also, this formula shows that each $S_{kj}$ belongs
 to a certain disk.

 For a fixed $k$, the operator $R_{kj}$ which appears in 
 \eqref{algoritm} is the row contraction
 $$R_{kj}=\left[
 \begin{array}{cccc}
 \Gamma _{k,k+1}, & D_{\Gamma ^*_{k,k+1}}\Gamma _{k,k+2},  
 & \ldots ,& D_{\Gamma ^*_{k,k+1}}\ldots D_{\Gamma ^*_{k,j-1}}\Gamma _{kj}
 \end{array}\right].$$
 Analogously, for a fixed $j$,
 the operator $C_{kj}$ is the column contraction
 $$C_{kj}=\left[
 \begin{array}{cccc}
 \Gamma _{j-1,j}, & \Gamma _{j-2,j}D_{\Gamma _{j-1,j}},  
 & \ldots ,& \Gamma _{kj}D_{\Gamma _{k-1,j}}\ldots D_{\Gamma _{j-1,j}}
 \end{array}\right]^t,$$
 where $"t"$ stands for matrix transpose.
 The operators $U_{ij}$ are defined by the recursion:
 $U_{kk}=1$ and for $j>k$,
 $$U_{kj}=U_j(\Gamma _{j,j+1})U_j(\Gamma _{j,j+2})\ldots
 U_j(\Gamma _{kj})(U_{k+1,j}\oplus I_{\cD _{\Gamma ^*_{kj}}}),$$
 where the subscript $j$ at $U(\Gamma _{j,j+l})$ means that for
 $1\leq l\leq j-k$ the unitary operator 
 $U_j(\Gamma _{k,k+l})$ is defined from 
 $$\left(\oplus _{m=1}^{l-1}\cD _{\Gamma _{k+1,k+m}}
 \right)\oplus (\cD _{\Gamma _{k+1,k+l}}
 \oplus \cD _{\Gamma ^*_{k,k+l}})\oplus
 \left(\oplus _{m=l+1}^j\cD _{\Gamma _{k,k+m}}\right)$$
 into 
 $$\left(\oplus _{m=1}^{l-1}\cD _{\Gamma _{k+1,k+m}}
 \right)\oplus (\cD _{\Gamma ^*_{k,k+l-1}}
 \oplus \cD _{\Gamma _{k,k+l}})\oplus
 \left(\oplus _{m=l+1}^j\cD _{\Gamma _{k,k+m}}\right)$$
 by the formula
 $$U_j(\Gamma _{k,k+l})=I\oplus U(\Gamma _{k,k+l})\oplus I.$$
 We note that the above formula for $U_{kj}$ comes from the 
 familiar Euler factorization of $SO(N)$, \cite{Mu}. 

 We obtain the following result.

 \begin{theorem}\label{param}
 There exists a one-to-one correspondence between the set
 of linear completely positive maps on $\cM_n$ and the 
 families $\Gamma =\{\Gamma _{kj}\mid 1\leq k\leq j\leq n^2\}$
 of complex numbers such that  
 $\Gamma _{kk}\geq 0$ for $1\leq k\leq n^2,$
 and for $1\leq k<j\leq n^2$,
 $\Gamma _{kj}\in \cL(\cD_{\Gamma _{k+1,j}},\cD_{\Gamma ^*_{k,j-1}})$
 are contractions. The correspondence is 
 given by \eqref{filtru} and \eqref{algoritm}.
 \end{theorem}

 This result can be rephrased as a Schur-Cohn type test for
 complete positivity.

 \begin{algorithm}\label{SC}
Consider a linear map
$\Phi $ on $\cM_n$. The complete positivity of 
$\Phi $ can be verified as follows:

$(1)$ Consider the matrix $S=S_{\Phi }$
given by formula \eqref{baza}. 

$(2)$ Check $\Gamma _{kk}\geq 0$ for each $k$. If for some $k$, 
$\Gamma _{kk}<0$, then 
$\Phi $ is not completely positive. If for some $k$, $\Gamma _{kk}=0$,
then the whole $k$th row (and column) of $S$ must be zero.

$(3)$ Calculate the numbers $\Gamma _{kj}$ according to  
formula \eqref{algoritm}. At each step check
the condition $|\Gamma _{kj}|\leq 1$ and keep track of the 
compatibility condition 
$\Gamma _{kj}\in \cL(\cD_{\Gamma _{k+1,j}},\cD_{\Gamma ^*_{k,j-1}})$.
If this can be done for all indices $kj$, then $\Phi $ is completely
positive. Otherwise, $\Phi $ is not completely positive.
\end{algorithm}

We illustrate the applicability of this algorithm for the 
case of completely positive maps on $\cM_2$.

 \begin{example}\label{m2}
 {\rm
 A detailed analysis of quantum binary channels is given in 
 \cite{RSW}. We show here how Theorem~\ref{param} relates to that 
 analysis. It is showed in \cite{KR} that any 
 quantum binary channel $\Phi $ has a representation
 $$\Phi (A)=U[\Phi _{{\bf t}, {\bf \Lambda }}(VAV^*)]U^*,$$
 where $U,V \in U(2)$ and 
 $\Phi _{{\bf t}, {\bf \Lambda }}$
 has the matrix representation
 $${\bf T}=\left[\begin{array}{cccc}
 1 & 0 & 0 & 0 \\
 t_1 &  \lambda _1 & 0 & 0 \\
 t_2 & 0 & \lambda _2 & 0 \\
 t_3 & 0 & 0 & \lambda _3 
 \end{array}\right]$$
 with respect to the Pauli basis 
 $\{I,\sigma _x,\sigma _y,\sigma _z\}$ of $\cM_2$.
 We can obtain (formula (26) in \cite{RSW}) that
 $$S_{\Phi _{{\bf t}, {\bf \Lambda }}}=
 \frac{1}{2}\left[\begin{array}{cccc}
 1+t_3+\lambda _3 & t_1-it_2 & 0 & \lambda _1+\lambda _2 \\
 t_1+it_2 & 1-t_3-\lambda _3 & \lambda _1-\lambda _2 & 0 \\
 0 & \lambda _1-\lambda _2 & 1+t_3-\lambda _3 & t_1-it_2 \\
 \lambda _1+\lambda _2 & 0 & t_1+it_2 & 1-t_3+\lambda _3
 \end{array}\right].
 $$
 Similarly, by formula (27) in \cite{RSW}, 
 $$S_{\hat{\Phi }_{{\bf t}, {\bf \Lambda }}}=
 \frac{1}{2}\left[\begin{array}{cccc}
 1+t_3+\lambda _3 & 0 & t_1+it_2 & \lambda _1+\lambda _2 \\
 0 & 1+t_3-\lambda _3 & \lambda _1-\lambda _2 & t_1+it_2 \\
 t_1-it_2 & \lambda _1-\lambda _2 & 1-t_3-\lambda _3 & 0 \\
 \lambda _1+\lambda _2 & t_1+it_2 & 0 & 1-t_3+\lambda _3
 \end{array}\right].
 $$
 It is slightly more convenient to deal with 
 $S=[S_{kj}]_{k,j=1}^4=2S_{\hat{\Phi }_{{\bf t}, {\bf \Lambda }}}$.
 Formula \eqref{algoritm} gives:
 $$S_{11}=\Gamma _{11}=1+t_3+\lambda _3;\quad \quad 
 S_{22}=\Gamma _{22}=1+t_3-\lambda _3;$$
 $$S_{33}=\Gamma _{33}=1-t_3-\lambda _3;\quad \quad 
 S_{44}=\Gamma _{44}=1-t_3+\lambda _3;$$
 $$\Gamma _{12}=0, \quad \Gamma _{34}=0;$$
 $$S_{23}=\Gamma _{22}^{1/2}\Gamma _{23}\Gamma _{33}^{1/2},$$
 so that 
 $$\Gamma _{23}=\frac{\lambda _1-\lambda _2}{(1+t_3-\lambda _3)^{1/2}
 (1-t_3-\lambda _3)^{1/2}};$$
 $$S_{13}=\Gamma ^{1/2}_{11}\Gamma _{13}D_{\Gamma _{23}}
 \Gamma ^{1/2}_{33},$$
 so that 
 $$\Gamma _{13}=\frac{(t_1+it_2)(1+t_3-\lambda _3)^{1/2}}
 {((1+t_3-\lambda _3)
 (1-t_3-\lambda _3)-(\lambda _1-\lambda _2)^2)^{1/2}
 (1+t_3+\lambda _3)^{1/2}};$$
 $$S_{24}=\Gamma ^{1/2}_{22}D_{\Gamma ^*_{23}}\Gamma _{24}
 \Gamma ^{1/2}_{44},$$
 so that 
 $$\Gamma _{24}=\frac{(t_1+it_2)(1-t_3-\lambda _3)^{1/2}}
 {((1+t_3-\lambda _3)
 (1-t_3-\lambda _3)-(\lambda _1-\lambda _2)^2)^{1/2}
 (1-t_3+\lambda _3)^{1/2}}.$$
 Finally, 
 $$S_{14}=\Gamma ^{1/2}_{11}(-\Gamma _{13}\Gamma ^*_{23}\Gamma _{24}+
 D_{\Gamma ^*_{13}}\Gamma _{14}D_{\Gamma _{24}})\Gamma ^{1/2}_{44}.$$
 By now, the formula for $\Gamma _{14}$ becomes quite intricate, but there
 in no problem to write it explicitly. We deduce that 
 ${\Phi }_{{\bf t}, {\bf \Lambda }}$ is completely positive 
 if and only if the following eigth inequalities hold:
 $$\Gamma _{kk}\geq 0, \quad k=1,\ldots ,4,$$
 $$|\Gamma _{23}|\leq 1,\quad |\Gamma _{13}|\leq 1,\quad
 |\Gamma _{24}|\leq 1,\quad |\Gamma _{14}|\leq 1.$$
 Also, we know what happens in the degenerate cases.
 Thus, the implication of $\Gamma _{kk}=0$ for some $k$ 
 on the structure of ${\Phi }_{{\bf t}, {\bf \Lambda }}$ is clear.
 Also, if $|\Gamma _{23}|=1$, 
 then necessarily $t_1=t_2=0$ and 
 $\lambda _1+\lambda _2=(1+t_3+\lambda _3)^{1/2}\Gamma _{14}
 (1-t_3+\lambda _3)^{1/2}$ for some contraction $\Gamma _{14}$.
 If either $|\Gamma _{13}|=1$ or $|\Gamma _{24}|=1$, then 
 necessarily 
$\Gamma _{14}=0$ and 
 $S_{14}=\Gamma ^{1/2}_{11}
(-\Gamma _{13}\Gamma ^*_{23}\Gamma _{24})\Gamma ^{1/2}_{44}$.

 We notice that this result is of about the same nature as that 
 in \cite{RSW}. This is because the first step 
 of \eqref{algoritm} is precisely Lemma~6 in \cite{RSW} which is 
 used for the analysis in \cite{RSW}. If we used the block version 
 of \eqref{algoritm} then we would deduce precisely Theorem~1
 of \cite{RSW}. What we basically have done here is that we used
 \eqref{algoritm} in order to deduce in a systematic way 
 the condition that $R_{{\Phi }_{{\bf t}, {\bf \Lambda }}}$
 in Theorem~1 of \cite{RSW} is a contraction. One advantage of 
 doing this is that it works in higher dimensions.

We also have to note that the correspondence between $S_{\Phi }$ and 
the parameters $\Gamma $ si nonlinear. Only for the first
step the correspondence is affine and therefore can be used in the 
analysis of extreme points in the case $n=2$, as it was done in \cite{RSW}. 
This seems to be 
unclear for $n\geq 2$. 
 }\qed
 \end{example}

We conclude with the presentation of so-called lattice structures that can be 
associated to completely positive maps on $\cM_n$. 
This comes from the remark that $S_{\Phi }$ has displacement 
structure as described in \cite{SCK} and the general 
lattice structures associated to matrices
with displacement structure in \cite{SCK} can be used in our
particular case. We can omit the details. In Figure~1 we show the
lattice structure of completely positive maps on $\cM_2$.

\begin{figure}[h]
\setlength{\unitlength}{2000sp}%
\begingroup\makeatletter\ifx\SetFigFont\undefined%
\gdef\SetFigFont#1#2#3#4#5{%
  \reset@font\fontsize{#1}{#2pt}%
  \fontfamily{#3}\fontseries{#4}\fontshape{#5}%
  \selectfont}%
\fi\endgroup%
\begin{picture}(7926,3912)(187,-3211)
\thinlines
\put(450, 409){\makebox(0,0)[lb]{\smash{\SetFigFont{12}{14.4}{\rmdefault}{\mddefault}{\updefault}\mbox{\small $\Gamma ^{1/2}_{44}$}}}}
\put(2350, 409){\makebox(0,0)[lb]{\smash{\SetFigFont{12}{14.4}{\rmdefault}{\mddefault}{\updefault}\mbox{\small $\Gamma ^{1/2}_{33}$}}}}
\put(5350, 409){\makebox(0,0)[lb]{\smash{\SetFigFont{12}{14.4}{\rmdefault}{\mddefault}{\updefault}\mbox{\small $\Gamma ^{1/2}_{22}$}}}}
\put(7300, 409){\makebox(0,0)[lb]{\smash{\SetFigFont{12}{14.4}{\rmdefault}{\mddefault}{\updefault}\mbox{\small $\Gamma ^{1/2}_{11}$}}}}

\put(301, 89){\circle{212}}
\put(6301, 89){\circle{212}}
\put(3301, 89){\circle{212}}
\put(601,239){\line( 0, 1){  0}}
\put(601,239){\line( 1, 0){1200}}
\put(1801,239){\line( 0,-1){1200}}
\put(1801,-961){\line(-1, 0){1200}}
\put(601,-961){\line( 0, 1){1200}}
\put(2101,-1861){\framebox(1200,1200){}}
\put(3601,-961){\framebox(1200,1200){}}
\put(5101,-1861){\framebox(1200,1200){}}
\put(6601,239){\framebox(0,0){}}
\put(6601,-961){\framebox(1200,1200){}}
\put(3601,-2761){\framebox(1200,1200){}}
\put(301, 89){\line( 1, 0){7800}}
\put(901,-811){\line( 1, 0){6600}}
\put(2401,-1711){\line( 1, 0){3600}}
\put(3901,-2611){\line( 1, 0){600}}
\put(901, 89){\vector( 2,-3){600}}
\put(901,-811){\vector( 2, 3){600}}
\put(3901, 89){\vector( 2,-3){600}}
\put(3901,-811){\vector( 2, 3){600}}
\put(6901, 89){\vector( 2,-3){600}}
\put(6901,-811){\vector( 2, 3){600}}
\put(2401,-811){\vector( 2,-3){600}}
\put(2401,-1711){\vector( 2, 3){600}}
\put(3901,-1711){\vector( 2,-3){600}}
\put(3901,-2611){\vector( 2, 3){600}}
\put(5401,-1711){\vector( 2, 3){600}}
\put(5401,-811){\vector( 2,-3){600}}
\put(5101, 89){\vector( 0, 1){600}}
\put(8101, 89){\vector( 0, 1){600}}
\put(2101, 89){\vector( 0, 1){600}}
\put(301, 14){\line( 0, 1){150}}
\put(226, 89){\line( 1, 0){ 75}}
\put(3301, 14){\line( 0, 1){150}}
\put(3301,164){\line( 0, 1){ 75}}
\put(6301, 14){\line( 0, 1){150}}
\put(3001, 89){\vector( 1, 0){225}}
\put(6001, 89){\vector( 1, 0){225}}
\put(301,689){\vector( 0,-1){450}}
\put(3301,689){\vector( 0,-1){450}}
\put(6301,689){\vector( 0,-1){450}}
\put(901, 89){\vector( 1, 0){600}}
\put(901,-811){\vector( 1, 0){525}}
\put(2401,-811){\vector( 1, 0){525}}
\put(2401,-1711){\vector( 1, 0){525}}
\put(3901, 89){\vector( 1, 0){525}}
\put(3901,-811){\vector( 1, 0){600}}
\put(3901,-1711){\vector( 1, 0){600}}
\put(3901,-2611){\vector( 1, 0){525}}
\put(5401,-811){\vector( 1, 0){600}}
\put(5401,-1711){\vector( 1, 0){600}}
\put(6901, 89){\vector( 1, 0){600}}
\put(6901,-811){\vector( 1, 0){600}}
\put(901,-1411){\makebox(0,0)[lb]{\smash{\SetFigFont{12}{14.4}{\rmdefault}{\mddefault}{\updefault}\mbox{\small $U(\Gamma _{34})$}}}}
\put(2401,-2311){\makebox(0,0)[lb]{\smash{\SetFigFont{12}{14.4}{\rmdefault}{\mddefault}{\updefault}\mbox{\small $U(\Gamma _{24})$}}}}
\put(3901,-1411){\makebox(0,0)[lb]{\smash{\SetFigFont{12}{14.4}{\rmdefault}{\mddefault}{\updefault}\mbox{\small $U(\Gamma _{23})$}}}}
\put(3901,-3211){\makebox(0,0)[lb]{\smash{\SetFigFont{12}{14.4}{\rmdefault}{\mddefault}{\updefault}\mbox{\small $U(\Gamma _{14})$}}}}
\put(5401,-2311){\makebox(0,0)[lb]{\smash{\SetFigFont{12}{14.4}{\rmdefault}{\mddefault}{\updefault}\mbox{\small $U(\Gamma _{13})$}}}}
\put(6901,-1411){\makebox(0,0)[lb]{\smash{\SetFigFont{12}{14.4}{\rmdefault}{\mddefault}{\updefault}\mbox{\small $U(\Gamma _{12})$}}}}
\end{picture}
\caption{\small Lattice structure for completely positive maps on $\cM_2$}
\end{figure}
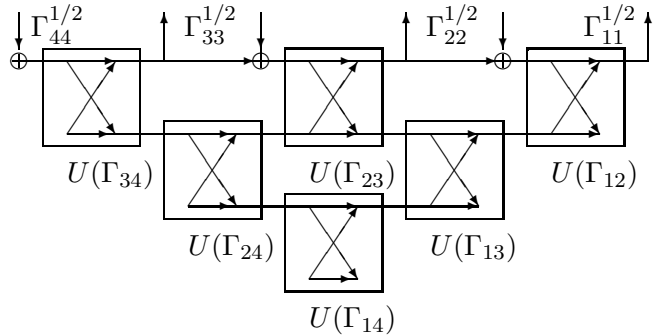

\end{document}